\documentclass[Svgc]{Svgc}
\usepackage{graphicx,psfrag,subfigure}
\usepackage{latexsym}
\usepackage{cite}
\usepackage{mathrsfs,amsmath,amsfonts,amssymb}

\usepackage{cases}
\usepackage{dsfont}
\usepackage{curves}
\usepackage{pifont}
\usepackage{array}
\usepackage{booktabs}
\usepackage[all,pdf]{xy}

\journalname{}

\newtheorem{cla}[subsection]{Claim}

\newtheorem{cor}[theorem]{Corollary}
\newtheorem{rem}{Remark}

\begin{document}

\title{Component factors in $K_{1,r}$-free graphs
\thanks{This work is supported by the National Natural Science Foundation of China (Grant No. 11871280, 11971349 and U1811461), the Natural Science Foundation of Guangdong Province (Grant No.~2020B1515310009) and Qinglan Project of Jiangsu Province.}
}
\mail{Xiaoyan Zhang (e-mail: royxyzhang@gmail.com, zhangxiaoyan@njnu.edu.cn)}
\author{Guowei Dai\inst{1} \and Zanbo Zhang\inst{2} \and Xiaoyan Zhang \inst{1}}

\titlerunning{Component factors in $K_{1,r}$-free graphs}
\authorrunning{G. Dai, X. Zhang and Z. Zhang}

\institute{
School of Mathematical Science \& Institute of Mathematics,
Nanjing Normal University, Nanjing 210023, P.R. China.
\and
School of Statistics \& Mathematics, and Institute of Artificial Intelligence \& Deep Learning, Guangdong University of Finance \& Economics, Guangzhou 510630, P.R. China.
}
\maketitle
\begin{abstract}
A graph is said to be $K_{1,r}$-free if it does not contain an induced subgraph isomorphic to $K_{1,r}$.
An $\mathcal{F}$-factor is a spanning subgraph $H$ such that each connected component of $H$ is isomorphic to some graph in $\mathcal{F}$.
In particular, $H$ is called an $\{P_2,P_3\}$-factor of $G$ if
$\mathcal{F}=\{P_2,P_3\}$;
$H$ is called an $\mathcal{S}_n$-factor of $G$ if $\mathcal{F}=\{K_{1,1},K_{1,2},K_{1,3},...,K_{1,n}\}$, where $n\geq2$.
A spanning subgraph of a graph $G$ is called a $\mathcal{P}_{\geq k}$-factor of $G$
if its each component is isomorphic to a path of order at least $k$, where $k\geq2$.
A graph $G$ is called a $\mathcal{F}$-factor covered graph if there is a
$\mathcal{F}$-factor of $G$ including $e$ for any $e\in E(G)$.
In this paper, we give a minimum degree condition for a $K_{1,r}$-free graph to have an $\mathcal{S}_n$-factor and a $\mathcal{P}_{\geq 3}$-factor, respectively.
Further, we obtain sufficient conditions for $K_{1,r}$-free graphs to be
$\mathcal{P}_{\geq 2}$-factor, $\mathcal{P}_{\geq 3}$-factor or
$\{P_2,P_3\}$-factor covered graphs.
In addition, examples show that our results are sharp.
\end{abstract}

\section{Introduction}
We consider only finite simple graph, unless explicitly stated.
We refer to \cite{Bondy1982} for the notation and terminologies not defined here.
Let $G = (V(G), E(G))$ be a simple graph, where $V(G)$ and $E(G)$ denote the vertex set and the edge set of $G$, respectively.
A subgraph $H$ of $G$ is called a spanning subgraph of $G$ if $V(H)=V(G)$ and $E(H)\subseteq E(G)$.
A subgraph $H$ of $G$ is called an induced subgraph of $G$ if every pair of vertices in $H$ which are adjacent in $G$ are also adjacent in $H$.
For $v\in V(G)$, we use $d_{G}(v)$ and $N_{G}(v)$ to denote the degree of $v$ and the set of vertices adjacent to $v$ in $G$, respectively.
For $S\subseteq V(G)$, we write $N_{G}(S)=\cup_{v\in S}N_{G}(v)$.
We use $\delta(G)$ to denote the minimum degree of a graph $G$.
We denote by $K_m$ the complete graph of order $m$.

A graph of order $r+1$ is called $K_{1,r}$ with center
vertex $u$ if it has the vertex set $\{u,u_1 ,u_2 ,...,u_r\}$ and edge set
$\{uu_1,uu_2,...,uu_r\}$, where $r$ is an integer with $r\geq2$.
We use $S_n$ to denote the set $\{K_2,K_{1,2},K_{1,3},...,K_{1,n}\}$, where $n\geq2$.
For an integer $r\geq3$, we say that a graph $G$ is $K_{1,r}$-free if $G$ does not contain an induced subgraph isomorphic to $K_{1,r}$.
In particular, a graph is said to be claw-free if it is $K_{1,3}$-free.

For a family of connected graphs $\mathcal{F}$, a spanning subgraph $H$ of a graph $G$ is called an $\mathcal{F}$-factor of $G$ if each component of $H$ is isomorphic to some graph in $\mathcal{F}$.
In particular, $H$ is called an $\mathcal{S}_n$-factor of $G$ if $\mathcal{F}=\{K_2,K_{1,2},K_{1,3},...,K_{1,n}\}$.
A spanning subgraph $H$ of a graph $G$ is called a $\mathcal{P}_{\geq k}$-factor of $G$ if every component of $H$ is isomorphic to a path of order at least $k$.
For example, a $\mathcal{P}_{\geq 3}$-factor means a graph factor in which every component is a path of order at least three.
A graph $G$ is called a $\mathcal{P}_{\geq k}$-factor covered graph if there is a $\mathcal{P}_{\geq k}$-factor of $G$ covering $e$ for any $e\in E(G)$.

Since Tutte proposed the well known Tutte 1-factor theorem \cite{Tutte1952}, there are many results on graph factors
\cite{Akiyama1985,Amahashi1982,Egawa2018,Kaneko2001,Wang1994},
$\mathcal{S}_n$-factors \cite{Egawa2008,Kano2010,Kano2012,Yu1990},
and $\mathcal{P}_{\geq k}$-factors in claw-free graphs and cubic graphs
\cite{Ando2002,Kano2008,Kawarabayashi2002}.
More results on graph factors can be found in the survey papers and books in \cite{Akiyama1985,Plummer2007,Yu2009}.
We use $\omega(G)$, $i(G)$ to denote the number of components and isolated vertices of a graph $G$, respectively.
For $X\subseteq V (G)$, let $G[X]$ be the subgraph of $G$ induced by $X$ and define $G-X := G[V (G)-X]$.
Akiyama, Avis and Era \cite{Akiyama1980} proved the following theorem, which is a criterion for a graph to have a $\mathcal{P}_{\geq 2}$-factor.

\begin{theorem}\label{thm:1}(Akiyama, Avis and Era \cite{Akiyama1985})
A graph $G$ has a $\mathcal{P}_{\geq 2}$-factor if and only if $i(G-X)\leq 2|X|$ for all $X\subseteq V(G)$.
\end{theorem}

Amahashi, Kano \cite{Amahashi1982} and Las Vergnas \cite{Las1978} gave independently
two necessary and sufficient  conditions for a graph to have an $\mathcal{S}_n$-factor, which is a generalization of Theorem \ref{thm:1}.

\begin{theorem}\label{thm:5}(Amahashi and Kano \cite{Amahashi1982}; Las Vergnas \cite{Las1978})
Let $n$ be an integer with $n\geq2$. Then a graph $G$ has an $\mathcal{S}_n$-factor if and only if $i(G-X)\leq n|X|$ for all $X\subseteq V(G)$, or $|S|\leq n|N_{G}(S)|$ for every independent set $S$ of $G$.
\end{theorem}

In this paper, we investigate the existence of $\mathcal{S}_n$-factors in $K_{1,r}$-free graphs, and
obtain a minimum degree condition for a $K_{1,r}$-free graph to have an $\mathcal{S}_n$-factor as following.

\begin{theorem}\label{thm:1-1}
Let $r,n$ be two integers with $r\geq3$ and $n\geq2$, and let $G$ be a $K_{1,r}$-free graph.
If $\delta(G)\geq \left\lfloor\frac{r+n-2}{n}\right\rfloor$,
then $G$ has an $\mathcal{S}_n$-factor.
\end{theorem}

\begin{rem}
The minimum degree condition in Theorem \ref{thm:1-1} is best possible.
To see this, we construct a graph $G=(r-1)K_1\vee G_1$, where
$G_1\cong K_{\left\lfloor\frac{r+n-2}{n}\right\rfloor-1}$ and $r\geq3, n\geq2$ are two integers.
It is obvious that $\delta(G)=\left\lfloor\frac{r+n-2}{n}\right\rfloor-1$ and $G$ is $K_{1,r}$-free.
We choose $X=V(G_1)$, and so $|X|=\left\lfloor\frac{r+n-2}{n}\right\rfloor-1$.
Thus, we have
$$
i(G-X)=r-1>n\times\left(\left\lfloor\frac{r+n-2}{n}\right\rfloor-1\right)=n|X|.
$$
In view of Theorem \ref{thm:5}, $G$ has no $\mathcal{S}_n$-factor.
\end{rem}

Since the graph $P_k,~k\geq4$ itself has a $\{P_2,P_3\}$-factor, it is easy to verify that
a graph has a $\mathcal{P}_{\geq 2}$-factor if and only if it has a $\{K_2,K_{1,2}\}$-factor.
Let $n=2$, then by Theorem \ref{thm:1-1}, we can obtain the the following corollary immediately.

\begin{cor}\label{cor:1-1}
Let $r$ be an integer with $r\geq3$, and let $G$ be a $K_{1,r}$-free graph.
If $\delta(G)\geq \left\lfloor\frac{r}{2}\right\rfloor$,
then $G$ has a $\mathcal{P}_{\geq 2}$-factor.
\end{cor}

Kaneko \cite{Kaneko2003} introduced the concept of a sun and gave a characterization for a graph with a $\mathcal{P}_{\geq 3}$-factor.
It is perhaps the first characterization of graphs which have a path factor not including $P_2$.
Recently, Kano et al. \cite{Kano2004} presented a simpler proof for Kaneko's theorem \cite{Kaneko2003}.

A graph $H$ is called factor-critical if $H-\{v\}$ has a 1-factor for each $v\in V(H)$.
Let $H$ be a factor-critical graph and $V(H)=\{v_1,v_2,...,v_n\}$.
By adding new vertices $\{u_1,u_2,...,u_n\}$ together with new edges $\{v_{i}u_{i} : 1 \leq i\leq n\}$ to $H$, the resulting graph is called a sun.
Note that, according to Kaneko \cite{Kaneko2003}, we regard $K_1$ and $K_2$ also as a sun, respectively.
Usually, the suns other than $K_1$ are called big suns.
It is called a sun component of $G-X$ if the component of $G-X$ is isomorphic to a sun.
We denote by $sun(G-X)$ the number of sun components in $G-X$.

\begin{theorem}\label{thm:2}(Kaneko \cite{Kaneko2003})
A graph $G$ has a $\mathcal{P}_{\geq 3}$-factor if and only if $sun(G-X)\leq 2|X|$ for all $X\subseteq V(G)$.
\end{theorem}

From Theorem \ref{thm:2}, Kaneko \cite{Kaneko2003} immediately deduce the following corollary.

\begin{cor}\label{thm:6}(Kaneko \cite{Kaneko2003})
Every claw-free graph with minimum degree at least 2 has a $\mathcal{P}_{\geq 3}$-factor.
\end{cor}

We extend Corollary \ref{thm:6} and obtain a minimum degree condition for a $K_{1,r}$-free graph admitting a $\mathcal{P}_{\geq 3}$-factor.

\begin{theorem}\label{thm:1-2}
Let $r$ be an integer with $r\geq3$, and let $G$ be a $K_{1,r}$-free graph.
If $\delta(G)\geq \left\lfloor\frac{r}{2}\right\rfloor+1$,
then $G$ has a $\mathcal{P}_{\geq 3}$-factor.
\end{theorem}

\begin{rem}
The minimum degree condition in Theorem \ref{thm:1-2} is best possible.
To see this, we construct a graph $G=(r-1)K_2\vee K_{\left\lfloor\frac{r}{2}\right\rfloor-1}$, where $r\geq3$ is an integer.
It is obvious that $\delta(G)=\left\lfloor\frac{r}{2}\right\rfloor$ and $G$ is $K_{1,r}$-free.
We choose $X=V(K_{\left\lfloor\frac{r}{2}\right\rfloor-1})$, and so $|X|=\left\lfloor\frac{r}{2}\right\rfloor-1$.
Thus, we have
$$
sun(G-X)=r-1>2\times\left(\left\lfloor\frac{r}{2}\right\rfloor-1\right)=2|X|.
$$
In view of Theorem \ref{thm:2}, $G$ has no $\mathcal{P}_{\geq 3}$-factor.
\end{rem}

Zhang and Zhou \cite{Zhang2009} proposed the concept of path-factor covered graphs, which is a generalization of matching covered graphs.
They also obtained a characterization for $\mathcal{P}_{\geq 2}$-factor and $\mathcal{P}_{\geq 3}$-factor covered graphs, respectively.

\begin{theorem}(Zhang et al. \cite{Zhang2009})\label{thm:3}
Let $G$ be a connected graph. Then $G$ is a $\mathcal{P}_{\geq 2}$-factor covered graph if and only if $i(G-S)\leq 2|S|-\varepsilon_1(S)$ for all
$S\subseteq V(G)$, where $\varepsilon_1(S)$ is defined by
\begin{equation*}
\varepsilon_1(S)=
\begin{cases}
2\quad \ \ &if~S\neq\emptyset~and~S~is~not~an~independent~set;\\
1\quad \ \ &if~S~is~a~nonempty~independent~set~and~there~exists\\
~\quad \ \ &a~nontrivial~component~of~$G$-$S$;\\
0\quad \ \ &otherwise.
\end{cases}
\end{equation*}
\end{theorem}

\begin{theorem}(Zhang et al. \cite{Zhang2009})\label{thm:4}
Let $G$ be a connected graph. Then $G$ is a $\mathcal{P}_{\geq 3}$-factor covered graph if and only if $sun(G-S)\leq 2|S|-\varepsilon_2(S)$ for all
$S\subseteq V(G)$, where $\varepsilon_2(S)$ is defined by
\begin{equation*}
\varepsilon_2(S)=
\begin{cases}
2\quad \ \ &if~S\neq\emptyset~and~S~is~not~an~independent~set;\\
1\quad \ \ &if~S~is~a~nonempty~independent~set~and~there~exists~a\\
~\quad \ \ & non-sun~component~of~$G$-$S$;\\
0\quad \ \ &otherwise.
\end{cases}
\end{equation*}
\end{theorem}

The first author of this paper has verified that a claw-free graph $G$ is
a $\mathcal{P}_{\geq 2}$-factor covered graph if $\delta(G)\geq 2$ \cite{Dai}.
We extend the above result and obtain a minimum degree condition for a $K_{1,r}$-free graph being a $\mathcal{P}_{\geq 2}$-factor covered graph.

\begin{theorem}\label{thm:2-1}
Let $r$ be an integer with $r\geq3$, and let $G$ be a $K_{1,r}$-free graph.
If $\delta(G)\geq \left\lfloor\frac{r}{2}\right\rfloor+1$,
then $G$ is a $\mathcal{P}_{\geq 2}$-factor covered graph.
\end{theorem}

\begin{rem}
The minimum degree condition in Theorem \ref{thm:2-1} is best possible.
In order to show this, we construct a graph $G=(r-1)K_1\vee K_{\left\lfloor\frac{r}{2}\right\rfloor}$, where $r\geq3$ is an integer.
It is obvious that $\delta(G)=\left\lfloor\frac{r}{2}\right\rfloor$ and $G$ is $K_{1,r}$-free.
We choose $S=V(K_{\left\lfloor\frac{r}{2}\right\rfloor})$, and so $|S|=\left\lfloor\frac{r}{2}\right\rfloor$.
Note that $\varepsilon_1(S)=2$.
Thus, we have
$$
i(G-S)=r-1>2\times\left\lfloor\frac{r}{2}\right\rfloor-2=
2|S|-\varepsilon_1(S).
$$
In terms of Theorem \ref{thm:3}, $G$ is not a $\mathcal{P}_{\geq 2}$-factor covered graph.
\end{rem}

The first author of this paper has also verified that a claw-free graph $G$ is
a $\mathcal{P}_{\geq 3}$-factor covered graph if $\delta(G)\geq 3$ \cite{Dai}.
We extend the above result and obtain a minimum degree condition for a $K_{1,r}$-free graph being a $\mathcal{P}_{\geq 3}$-factor covered graph.

\begin{theorem}\label{thm:2-2}
Let $r$ be an integer with $r\geq3$, and let $G$ be a $K_{1,r}$-free graph.
If $\delta(G)\geq \left\lfloor\frac{r}{2}\right\rfloor+2$,
then $G$ is a $\mathcal{P}_{\geq 3}$-factor covered graph.
\end{theorem}

\begin{rem}
The minimum degree condition in Theorem \ref{thm:2-2} is best possible.
In order to show this, we construct a graph $G=(r-1)K_2\vee K_{\left\lfloor\frac{r}{2}\right\rfloor}$, where $r\geq3$ is an integer.
It is obvious that $\delta(G)=\left\lfloor\frac{r}{2}\right\rfloor+1$ and $G$ is $K_{1,r}$-free.
We choose $S=V(K_{\left\lfloor\frac{r}{2}\right\rfloor})$, and so $|S|=\left\lfloor\frac{r}{2}\right\rfloor$.
Note that $\varepsilon_2(S)=2$.
Thus, we have
$$
sun(G-S)=r-1>2\times\left\lfloor\frac{r}{2}\right\rfloor-2=
2|S|-\varepsilon_2(S).
$$
In terms of Theorem \ref{thm:4}, $G$ is not a $\mathcal{P}_{\geq 3}$-factor covered graph.
\end{rem}

Let $n=2$, then $\mathcal{S}_n=\{P_2,P_3\}$.
The following is a main result in \cite{Yu1991} and \cite{Yu1987}.

\begin{theorem}\label{thm:5}
Let $G$ be a connected graph. $G$ is a $\{P_2,P_3\}$-factor covered graph
if and only if $i(G-S)\leq 2|S|-\varepsilon_3(S)$ for all $S\subseteq V(G)$,
where $\varepsilon_3(S)$ is defined by
\begin{equation*}
\varepsilon_3(S)=
\begin{cases}
3\quad \ \ &if~nonempty~set~X~is~not~independent;\\
0\quad \ \ &otherwise.
\end{cases}
\end{equation*}
\end{theorem}

Note that a graph is a $\mathcal{P}_{\geq 2}$-factor covered graph not equivalent to a $\{P_2,P_3\}$-factor covered graph.
For example, $P:=x_1x_2x_3x_4$ is a $\mathcal{P}_{\geq 2}$-factor covered graph and it has no $\{P_2,P_3\}$-factor to cover edge $x_2x_3$.
We give a sufficient condition for a $K_{1,r}$-free graph to be a $\{P_2,P_3\}$-factor covered graph.

\begin{theorem}\label{thm:2-3}
Let $r$ be an integer with $r\geq3$, and let $G$ be a $K_{1,r}$-free graph.
If $\delta(G)\geq \left\lceil\frac{r}{2}\right\rceil+1$,
then $G$ is a $\{P_2,P_3\}$-factor covered graph.
\end{theorem}

\begin{rem}
The minimum degree condition in Theorem \ref{thm:2-3} is best possible.
In order to show this, we construct a graph
$G=(r-1)K_1\vee K_{\left\lceil\frac{r}{2}\right\rceil}$, where $r\geq3$ is an integer.
It is obvious that $\delta(G)=\left\lceil\frac{r}{2}\right\rceil$ and $G$ is $K_{1,r}$-free.
We choose $S=V(K_{\left\lceil\frac{r}{2}\right\rceil})$, and so $|S|=\left\lceil\frac{r}{2}\right\rceil$.
Note that $\varepsilon_3(S)=3$.
Thus, we have
$$
i(G-S)=r-1>2\times\left\lceil\frac{r}{2}\right\rceil-3=
2|S|-\varepsilon_3(S).
$$
In terms of Theorem \ref{thm:5}, $G$ is not a $\{P_2,P_3\}$-factor covered graph.
\end{rem}

\section{Proof of Theorem \ref{thm:1-1}}
By contradiction, suppose that Theorem \ref{thm:1-1} is not hold.
Then, there exists a $K_{1,r}$-free graph $G$ such that $\delta(G)\geq \left\lfloor\frac{r+n-2}{n}\right\rfloor$, and $G$ has no $\mathcal{S}_n$-factor.
It is easy to find that $G$ has a connected component $G'$ such that $G'$ has no $\mathcal{S}_n$-factor.
Note that $G'$ is also a $K_{1,r}$-free graph and $\delta(G')\geq \left\lfloor\frac{r+n-2}{n}\right\rfloor\geq1$.
By Theorem \ref{thm:5}, there exists a subset $X\subseteq V(G')$ such that $i(G'-X)>n|X|$.
In terms of the integrality of $i(G'-X)$, we obtain that
\begin{equation}\label{eqn1-5}
i(G'-X)\geq n|X|+1.
\end{equation}

\begin{cla}\label{cla:1}
$X\neq\emptyset$.
\end{cla}
\proof
Suppose $X=\emptyset$, then $i(G')=i(G'-X)\geq n|X|+1=1$.
Since $G'$ is a connected graph, we obtain that $G'$ is a isolated vertex, which contradicts the minimum degree of $G'$.
This completes the proof of Claim \ref{cla:1}. \qed

\vspace{3mm}

By Claim \ref{cla:1}, $X\neq\emptyset$.
Let $X=\{x_1,x_2,...,x_k\}$, then $|X|=k$.
Let $Y=\{y_1,y_2,\ldots,y_m\}$ be the set of isolated vertices of $G'-X$.
Then, by (\ref{eqn1-5}), we have that
\begin{equation}\label{eqn1-6}
i(G'-X)=m\geq nk+1.
\end{equation}

Next, we construct a bipartite subgraph $F$ of $G'$ such that
$V(F)=V(X)\cup V(Y)$ and $x_iy_j\in E(F)$ if and only if $x_iy_j\in E(G')$
for any $i\in [1,k],~j\in [1,m]$.
Since $\delta(G')\geq \left\lfloor\frac{r+n-2}{n}\right\rfloor$,
we have $d_{X}(y_j)\geq \left\lfloor\frac{r+n-2}{n}\right\rfloor$ for any $j\in [1,m]$.
By (\ref{eqn1-6}), we can derive that
\begin{equation*}
|X|\geq d_{X}(y_j)\geq \left\lfloor\frac{r+n-2}{n}\right\rfloor,
\end{equation*}
and thus
\begin{eqnarray*}
|E(F)|
&=& \sum_{j=1}^m d_F(y_j)
\\
&\geq&  m\times\left\lfloor\frac{r+n-2}{n}\right\rfloor
\\
&\geq& (nk+1)\times\frac{r-1}{n}
\\
&=& kr+\frac{r-1}{n}.
\end{eqnarray*}
It follows immediately that
\begin{equation}\label{eqn1-7}
\frac{|E(F)|}{|X|}\geq \frac{k(r-1)+\frac{r-1}{n}}{k}=r-1+\frac{r-1}{nk}>r-1.
\end{equation}
By pigeonhole principle and (\ref{eqn1-7}), there exists $x_{i}\in X$ such that $d_F(x_i)\geq r$.
Let $J=N_F(x_i)\cap Y$, then $G'[\{x_i\}\cup J]\cong K_{1,|J|}$.
Since $G'$ is $K_{1,r}$-free, $|J|<r$, which is a contradiction to that $|J|=d_F(x_i)\geq r$.

This completes the proof of Theorem \ref{thm:1-1}.

\section{Proof of Theorem \ref{thm:1-2}}
By contradiction, suppose that Theorem \ref{thm:1-2} is not hold.
Then, there exists a $K_{1,r}$-free graph $G$ such that $\delta(G)\geq \left\lfloor\frac{r}{2}\right\rfloor+1$, and $G$ has no $\mathcal{P}_{\geq3}$-factor.
It is easy to find that $G$ has a connected component $G'$ such that $G'$ has no $\mathcal{P}_{\geq3}$-factor.
Note that $G'$ is also a $K_{1,r}$-free graph and $\delta(G')\geq \left\lfloor\frac{r}{2}\right\rfloor+1\geq2$.
By Theorem \ref{thm:2}, there exists a subset $X\subseteq V(G')$ such that $sun(G'-X)>2|X|$.
In terms of the integrality of $sun(G'-X)$, we obtain that
\begin{equation}\label{eqn1-1}
sun(G'-X)\geq 2|X|+1.
\end{equation}

\begin{cla}\label{cla:2}
$X\neq\emptyset$.
\end{cla}
\proof
Suppose $X=\emptyset$, then $sun(G')=sun(G'-X)\geq2|X|+1=1$.
On the other hand, $sun(G')\leq\omega(G')=1$ since $G'$ is a connected graph.
So, we obtain that $G'$ is a big sun, which contradicts the minimum degree of $G'$.
This completes the proof of Claim \ref{cla:2}. \qed

\vspace{3mm}

By Claim \ref{cla:2}, $X\neq\emptyset$.
Let $X=\{x_1,x_2,...,x_k\}$, then $|X|=k$.
Let $\{C_1 ,C_2 ,...,C_m\}$ be the set of sun components of $G'-X$.
Then, by (\ref{eqn1-1}), we have that
\begin{equation}\label{eqn1-2}
sun(G'-X)=m\geq2k+1.
\end{equation}
For any $i\in [1,m]$, as $C_i$ is a sun component, there exists $y_i\in V(C_i)$ such that $d_{C_i}(y_i)\leq 1$.
Let $Y=\{y_1,y_2,\ldots,y_m\}$.
Then we construct a bipartite subgraph $F$ of $G'$ such that
$V(F)=V(X)\cup V(Y)$ and $x_iy_j\in E(F)$ if and only if $x_iy_j\in E(G')$
for any $i\in [1,k],~j\in [1,m]$.
Note that for any $j\in [1,m]$, we have $d_{C_j}(y_j)\leq 1$.
Thus, for any $j\in [1,m]$,
\begin{equation}\label{eqn1-3}
d_F(y_j)=d_{G'}(y_j)-d_{C_j}(y_j)\geq \left\lfloor\frac{r}{2}\right\rfloor.
\end{equation}
By (\ref{eqn1-2}) and (\ref{eqn1-3}), we can derive that
\begin{equation*}
|X|\geq d_F(y_j)\geq \left\lfloor\frac{r}{2}\right\rfloor,
\end{equation*}
and thus
\begin{eqnarray*}
|E(F)|
&=& \sum_{j=1}^m d_F(y_j)
\\
&\geq&  m\times\left\lfloor\frac{r}{2}\right\rfloor
\\
&\geq& (2k+1)\times\frac{r-1}{2}
\\
&=& k(r-1)+\frac{r-1}{2}.
\end{eqnarray*}
It follows immediately that
\begin{equation}\label{eqn1-4}
\frac{|E(F)|}{|X|}\geq \frac{k(r-1)+\frac{r-1}{2}}{k}=r-1+\frac{r-1}{2k}>r-1.
\end{equation}
By pigeonhole principle and (\ref{eqn1-4}), there exists $x_{i}\in X$ such that $d_F(x_i)\geq r$.
Let $I=N_F(x_i)\cap Y$, then $G'[\{x_i\}\cup I]\cong K_{1,|I|}$.
Since $G'$ is $K_{1,r}$-free, $|I|<r$, which is a contradiction to that $|I|=d_F(x_i)\geq r$.
This completes the proof of Theorem \ref{thm:1-2}.

\section{Proof of Theorem \ref{thm:2-1}}
By contradiction, suppose that Theorem \ref{thm:2-1} is not hold.
Then, there exists a $K_{1,r}$-free graph $G$ such that $\delta(G)\geq \left\lfloor\frac{r}{2}\right\rfloor+1$,
and $G$ is not a $\mathcal{P}_{\geq2}$-factor covered graph.
It is easy to find that $G$ has a connected component $G'$ such that $G'$ is not a $\mathcal{P}_{\geq2}$-factor covered graph.
Note that $G'$ is also a $K_{1,r}$-free graph and $\delta(G')\geq \left\lfloor\frac{r}{2}\right\rfloor+1\geq2$.
By Theorem \ref{thm:3}, there exists a subset $S\subseteq V(G')$ such that $i(G'-S)>2|S|-\varepsilon_1(S)$.
In terms of the integrality of $i(G'-S)$, we obtain that
$i(G'-S)\geq 2|S|-\varepsilon_1(S)+1$.

We will distinguish two cases below to show that $G$ is a $\mathcal{P}_{\geq2}$-factor covered graph. \\

{\bf Case~1}. $S=\emptyset$.

In this case, by the definition of $\varepsilon_1(S)$, we have $\varepsilon_1(S)=0$.
It follows easily that
$$i(G')=i(G'-S)\geq2|S|-\varepsilon_1(S)+1=1.$$
Since $G'$ is a connected graph, we obtain that $G'$ is a isolated vertex, which contradicts the minimum degree of $G'$.
This completes the proof of Case 1. \\

{\bf Case~2}. $S\neq\emptyset$.

Let $|S|=k$ and $S=\{s_1,s_2,...,s_k\}$.
By the definition of $\varepsilon_1(S)$, we have $\varepsilon_1(S)\leq 2$.
It follows easily that
$$i(G'-S)\geq2|S|-\varepsilon_1(S) + 1\geq2|S|-1.$$
Let $X=\{x_1,x_2,\ldots,x_m\}$ be the set of isolated vertices of $G'-S$,
then $i(G'-S) = m\geq2k-1$.
Next, we construct a bipartite subgraph $H$ of $G'$ such that
$V(H)=V(S)\cup V(X)$ and $s_ix_j\in E(H)$ if and only if $s_ix_j\in E(G')$
for any $i\in [1,k],~j\in [1,m]$.
Since $\delta(G')\geq \left\lfloor\frac{r}{2}\right\rfloor+1$,
we have $d_{S}(x_j)\geq \left\lfloor\frac{r}{2}\right\rfloor+1$ for any $j\in [1,m]$.
Then, we can derive that
\begin{equation*}
|S|\geq d_{S}(x_j)\geq \left\lfloor\frac{r}{2}\right\rfloor+1,
\end{equation*}
and thus
\begin{eqnarray*}
|E(H)|
&=& \sum_{j=1}^m d_H(x_j)
\\
&\geq&  m\left(\left\lfloor\frac{r}{2}\right\rfloor+1\right)
\\
&\geq& (2k-1)\times\frac{r+1}{2}
\\
&=& kr+\left(k-\frac{r+1}{2}\right)
\\
&\geq& kr.
\end{eqnarray*}
It follows immediately that
\begin{equation}\label{eqn2-3}
\frac{|E(H)|}{|S|}\geq \frac{kr}{|S|}=r.
\end{equation}
By pigeonhole principle and (\ref{eqn2-3}), there exists $s_{i}\in S$ such that $d_H(s_i)\geq r$.
Let $J=N_H(s_i)\cap X$, then $G'[\{s_i\}\cup J]\cong K_{1,|J|}$.
Since $G'$ is $K_{1,r}$-free, $|J|<r$, which is a contradiction to that $|J|=d_H(s_i)\geq r$.
This completes the proof of Case 2.

Combining Case 1 and Case 2, Theorem \ref{thm:2-1} is proved.

\section{Proof of Theorem \ref{thm:2-2}}
By contradiction, suppose that Theorem \ref{thm:2-2} is not hold.
Then, there exists a $K_{1,r}$-free graph $G$ such that $\delta(G)\geq \left\lfloor\frac{r}{2}\right\rfloor+2$,
and $G$ is not a $\mathcal{P}_{\geq3}$-factor covered graph.
It is easy to find that $G$ has a connected component $G'$ such that $G'$ is not a $\mathcal{P}_{\geq3}$-factor covered graph.
Note that $G'$ is also a $K_{1,r}$-free graph and $\delta(G')\geq \left\lfloor\frac{r}{2}\right\rfloor+2\geq3$.
By Theorem \ref{thm:4}, there exists a subset $S\subseteq V(G')$ such that $sun(G'-S)>2|S|-\varepsilon_2(S)$.
In terms of the integrality of $sun(G'-S)$, we obtain that
$sun(G'-S)\geq 2|S|-\varepsilon_2(S)+1$.

We will distinguish two cases below to show that $G$ is a $\mathcal{P}_{\geq3}$-factor covered graph. \\

{\bf Case~1}. $S=\emptyset$.

In this case, by the definition of $\varepsilon_2(S)$, we have $\varepsilon_2(S)=0$.
It follows easily that
$$sun(G')=sun(G'-S)\geq2|S|-\varepsilon_2(S)+1=1.$$
On the other hand, $sun(G')\leq\omega(G')=1$ since $G'$ is a connected graph.
Combining the results above, we obtain that $G'$ is a big sun, which contradicts the minimum degree of $G'$.
This completes the proof of Case 1. \\

{\bf Case~2}. $S\neq\emptyset$.

Let $|S|=k$ and $S=\{s_1,s_2,...,s_k\}$.
By the definition of $\varepsilon_2(S)$, we have $\varepsilon_2(S)\leq 2$.
It follows easily that
$$sun(G'-S)\geq2|S|-\varepsilon_2(S) + 1\geq2|S|-1.$$
Let $\{C_1 ,C_2 ,...,C_m\}$ be the set of sun components of $G'-S$,
then $sun(G'-S) = m\geq2k-1$.
For any $i\in [1,m]$, as $C_i$ is a sun component, there exists $x_i\in V(C_i)$ such that $d_{C_i}(x_i)\leq 1$.
Let $X=\{x_1,x_2,\ldots,x_m\}$.
Then we construct a bipartite subgraph $H$ of $G'$ such that
$V(H)=V(S)\cup V(X)$ and $s_ix_j\in E(H)$ if and only if $s_ix_j\in E(G')$
for any $i\in [1,k],~j\in [1,m]$.
Note that for any $j\in [1,m]$, we have $d_{C_j}(x_j)\leq 1$.
Thus, for any $j\in [1,m]$,
\begin{equation}\label{eqn2-1}
d_H(x_j)=d_{G'}(x_j)-d_{C_j}(x_j)\geq \left\lfloor\frac{r}{2}\right\rfloor+1.
\end{equation}
By (\ref{eqn2-1}), we can derive that
\begin{equation*}
|S|\geq d_H(x_j)\geq \left\lfloor\frac{r}{2}\right\rfloor+1,
\end{equation*}
and thus
\begin{eqnarray*}
|E(H)|
&=& \sum_{j=1}^m d_H(x_j)
\\
&\geq&  m\left(\left\lfloor\frac{r}{2}\right\rfloor+1\right)
\\
&\geq& (2k-1)\times\frac{r+1}{2}
\\
&=& kr+\left(k-\frac{r+1}{2}\right)
\\
&\geq& kr.
\end{eqnarray*}
It follows immediately that
\begin{equation}\label{eqn2-2}
\frac{|E(H)|}{|S|}\geq \frac{kr}{|S|}=r.
\end{equation}
By pigeonhole principle and (\ref{eqn2-2}), there exists $s_{i}\in S$ such that $d_H(s_i)\geq r$.
Let $I=N_H(s_i)\cap X$, then $G'[\{s_i\}\cup I]\cong K_{1,|I|}$.
Since $G'$ is $K_{1,r}$-free, $|I|<r$, which is a contradiction to that $|I|=d_H(s_i)\geq r$.
This completes the proof of Case 2.

Combining Case 1 and Case 2, Theorem \ref{thm:2-2} is proved.

\section{Proof of Theorem \ref{thm:2-3}}
By contradiction, suppose that Theorem \ref{thm:2-3} is not hold.
Then, there exists a $K_{1,r}$-free graph $G$ such that $\delta(G)\geq \left\lceil\frac{r}{2}\right\rceil+1$,
and $G$ is not a $\{P_2,P_3\}$-factor covered graph.
It is easy to find that $G$ has a connected component $G'$ such that $G'$ is not a $\{P_2,P_3\}$-factor covered graph.
Note that $G'$ is also a $K_{1,r}$-free graph and $\delta(G')\geq \left\lceil\frac{r}{2}\right\rceil+1\geq3$.
By Theorem \ref{thm:5}, there exists a subset $S\subseteq V(G')$ such that $i(G'-S)>2|S|-\varepsilon_3(S)$.
In terms of the integrality of $i(G'-S)$, we obtain that
$i(G'-S)\geq 2|S|-\varepsilon_3(S)+1$.

We will distinguish two cases below to show that $G$ is a $\{P_2,P_3\}$-factor covered graph. \\

{\bf Case~1}. $|S|\leq1$.

In this case, by the definition of $\varepsilon_3(S)$, we have $\varepsilon_3(S)=0$.
It follows easily that
$$i(G')=i(G'-S)\geq2|S|-\varepsilon_3(S)+1=1.$$
Since $G'$ is a connected graph, we obtain that $G'$ is a isolated vertex, which contradicts the minimum degree of $G'$.
This completes the proof of Case 1. \\

{\bf Case~2}. $|S|\geq2$.

Let $|S|=k$ and $S=\{s_1,s_2,...,s_k\}$.
By the definition of $\varepsilon_3(S)$, we have $\varepsilon_3(S)\leq 3$.
It follows easily that
$$i(G'-S)\geq2|S|-\varepsilon_3(S)+1\geq2|S|-2.$$
Let $X=\{x_1,x_2,\ldots,x_m\}$ be the set of isolated vertices of $G'-S$,
then $i(G'-S) = m\geq2k-2$.
Next, we construct a bipartite subgraph $H$ of $G'$ such that
$V(H)=V(S)\cup V(X)$ and $s_ix_j\in E(H)$ if and only if $s_ix_j\in E(G')$
for any $i\in [1,k],~j\in [1,m]$.
Since $\delta(G')\geq \left\lceil\frac{r}{2}\right\rceil+1$,
we have $d_{S}(x_j)\geq \left\lceil\frac{r}{2}\right\rceil+1$ for any $j\in [1,m]$.
Then, we can derive that
\begin{equation*}
|S|\geq d_{S}(x_j)\geq \left\lceil\frac{r}{2}\right\rceil+1,
\end{equation*}
and thus
\begin{eqnarray*}
|E(H)|
&=& \sum_{j=1}^m d_H(x_j)
\\
&\geq&  m\left(\left\lceil\frac{r}{2}\right\rceil+1\right)
\\
&\geq& (2k-2)\times\left(\frac{r}{2}+1\right)
\\
&=& kr+2\times\left(k-\frac{r}{2}-1\right)
\\
&\geq& kr.
\end{eqnarray*}
It follows immediately that
\begin{equation}\label{eqn2-4}
\frac{|E(H)|}{|S|}\geq \frac{kr}{|S|}=r.
\end{equation}
By pigeonhole principle and (\ref{eqn2-4}), there exists $s_{i}\in S$ such that $d_H(s_i)\geq r$.
Let $J=N_H(s_i)\cap X$, then $G'[\{s_i\}\cup J]\cong K_{1,|J|}$.
Since $G'$ is $K_{1,r}$-free, $|J|<r$, which is a contradiction to that $|J|=d_H(s_i)\geq r$.
This contradiction completes the proof of Theorem \ref{thm:2-3}.

\begin{acknowledgement}
We are very grateful for the helpful comments from Professor Zhiquan Hu, which greatly help to improve the manuscript.
\end{acknowledgement}

\end{document}